\documentclass[oneside,reqno]{amsart}
\usepackage{amsfonts,amsmath,latexsym,verbatim,amscd,mathrsfs,color,array}
\usepackage[colorlinks=true]{hyperref}

\usepackage{amsmath,amssymb,amsthm,amsfonts,graphicx,color}
\usepackage{amssymb}
\usepackage{pdfsync}
\usepackage{epstopdf}

\newcommand{\ass}{\quad\mbox{as}\quad}

\newcommand{\inn}{{\quad\hbox{in } }}
\newcommand{\onn}{{\quad\hbox{on } }}
\newcommand{\ttt}{\tilde }

\newcommand{\TT}{{\mathcal T}  }

\newcommand{\nn}{ {\nabla}  }

\newcommand{\pp}{ {\partial} }

\newcommand{\vp}{\varphi}

\newcommand{\R} {\mathbb R}

\newcommand{\cuad}{{\sqcap\kern-.68em\sqcup}}

\newcommand{\DD}{{\mathcal D}}

\newcommand{\foral}{\quad\mbox{for all}\quad}
\newcommand{\ve}{\varepsilon}

\newcommand{\be}{\begin{equation}}
\newcommand{\ee}{\end{equation}}

\newcommand{\la}{\lambda}

\newcommand{\equ}[1]{(\ref{#1})}

\newcommand{\FF}{{\mathcal  F}  }

\newtheorem{lemma}{Lemma}[section]

\newtheorem{theorem}{Theorem}

\newtheorem{remark}{Remark}[section]
\newcommand{\bremark}{\begin{remark} \em}
\newcommand{\eremark}{\end{remark} }

\numberwithin{equation}{section}

\begin{document}

\title[Type II blow up for the 5-dimensional critical heat equation]{ Type II blow-up in the 5-dimensional energy critical heat equation}

\author[M. del Pino]{Manuel del Pino}
\address{\noindent   Department of Mathematical Sciences University of Bath,
Bath BA2 7AY, United Kingdom \\
and  Departamento de
Ingenier\'{\i}a  Matem\'atica-CMM   Universidad de Chile,
Santiago 837-0456, Chile}
\email{m.delpino@bath.ac.uk}

\author[M. del Pino]{Monica Musso}
\address{\noindent   Department of Mathematical Sciences University of Bath,
Bath BA2 7AY, United Kingdom \\
and Departamento de Matem\'aticas, Universidad Cat\'olica de Chile, Macul 782-0436, Chile}
\email{m.musso@bath.ac.uk}

\author[J. Wei]{Juncheng Wei}
\address{\noindent  Department of Mathematics University of British Columbia, Vancouver, BC V6T 1Z2, Canada
}  \email{jcwei@math.ubc.ca}

\begin{abstract}
We consider the Cauchy problem for the energy critical heat equation
\be\left\{
\begin{aligned}
u_t   & = \Delta u + |u|^{\frac 4{n-2}}u   \inn\ \R^n \times (0, T),  \\
u(\cdot,0) &  =u_0 \inn \R^n
\end{aligned}\right.
\label{01}\ee
in dimension $n=5$. More precisely we find that for given points $q_1, q_2,\ldots, q_k$ and any sufficiently small $T>0$ there is an initial condition $u_0$ such that the solution  $u(x,t)$  of \equ{01} blows-up at exactly those $k$ points with rates type II, namely
 with absolute size $ \sim (T-t)^{-\alpha} $ for $\alpha > \frac 34  $.
 The blow-up profile around each point
is of bubbling type, in the form of sharply scaled Aubin-Talenti bubbles.

\end{abstract}

\maketitle


\section{Introduction}
Many studies have been devoted to the analysis of blow-up phenomena in a semilinear heat equation
of the form
\be\left\{
\begin{aligned}
u_t   & = \Delta u+ |u|^{p-1}u   \inn\ \Omega \times (0, T),  \\
u &  =0 \onn \pp\Omega\times (0,T),\\
u(\cdot,0) &  =u_0 \inn \Omega.
\end{aligned}\right.
\label{F}\ee
where $p>1$,  and $\Omega$ is a bounded smooth domain in $\R^n$ (or entire space), starting with the seminal work by Fujita \cite{fujita} in the 1960's. A smooth solution of \equ{F} {\em blows-up at time $T$} if
$$ \lim_{t\to T} \| u(\cdot, t)\|_{L^\infty (\R^n) }    = + \infty . $$
We observe that for functions independent of the space variable the equation reduces to the ODE
$ u_t =  |u|^{p-1}u $, which is solved
for a suitable constant $c_p$ by the function $ u(t)= c_p (T-t)^{-\frac 1{p-1}} $ and it blows-up at time $T$. It is commonly said that the blow-up of a solution $u(x,t)$ is of type I if it happens at most at the ODE rate:
$$
\limsup_{t\to T} (T-t)^{\frac 1{p-1}} \| u(\cdot, t)\|_{L^\infty (\R^n) }    <  + \infty
$$
while the blow-up is said of type II if
$$
\limsup_{t\to T} (T-t)^{\frac 1{p-1}} \| u(\cdot, t)\|_{L^\infty (\R^n) }    =  + \infty
$$
Many results have predicted that type I is the ``typical'' or ``generic'' way in which blows-up  takes place for solutions of  equation \equ{F}. For instance it is known after a series of works, including \cite{gk1,gk2,gk3}, that type I is the only way possible if $p<p_S$ where $p_S$ is the critical Sobolev exponent,
$$p_S :=\begin{cases} \frac{n+2}{n-2}  & \hbox{ if } n\ge 3 \\ +\infty  &\hbox{ if } n=1,2.    \end{cases}
$$
Stability and genericity of type blow-up have been considered for instance in \cite{collot2,mz,mm2}.
Solutions with type II blow-up are in fact much harder to detect. The first example was discovered in \cite{hv,hv1}, for $p>p_{JL}$ where    $p_{JL}$ is the Joseph-Lundgreen exponent \cite{jl},
$$p_{JL}=\begin{cases} 1+{4\over n-4-2\,\sqrt{n-1}} & \text{if $n\ge11$}\\ + \infty, & \text{if $ n\le 10$}.\\  \end{cases} $$
See the book \cite{qs} for a survey of related results.
In fact, no type II blow-up is present for radial solutions if $p_S< p< p_{JL}$, while for radial positive solutions this is not possible
if $p= \frac{n+2}{n-2}$ \cite{fhv1}. Examples of nonradial positive blow-up solutions for $p> p_{JL}$  have been found in
\cite{collot4,collot1}.
 Formally a sign-changing solution with type II blow-up was predicted for $p=p_S$ in \cite{fhv1} and a rigorous radial example in \cite{schweyer} for $n=4$.  In this paper we shall exhibit a first example of type II blow-up for dimension $n=5$ and critical case    $p=p_S$ and  $p=p_S = \frac 73 $. Our construction does not depend on any symmetries, in fact it can be made in such a way that
 blow-up takes place simultaneously at any prescribed points  of space with a {\em bubbling profile}. We recall that all
 positive entire solutions of the equation
 $$
 \Delta u  +  |u|^{\frac 4{n-2}}u = 0 \inn \R^n
 $$
 are given by the family of {\em Aubin-Talenti bubbles}
 \be
 U_{\la,\xi} (x)  =   \la^{-\frac {n-2}2} U\left (\frac {x-\xi}{\la}   \right )
 \label{bubbles} \ee
 where $U(y)$ is the {\em standard bubble}
 $$
 U(y)  =   \alpha_n \left ( \frac 1{1+|y|^2}\right)^{\frac{n-2}2} , \quad \alpha_n = (n(n-2))^{\frac 1{n-2}}.
 $$
 The solutions we construct do change sign, and look at main order near the blow-up points as one of the bubbles \equ{bubbles} with time dependent parameters and $\mu(t)\to 0$
 as $t\to T$.
 \medskip
 Thus we consider the equation
\be\left\{
\begin{aligned}
u_t   & = \Delta u + |u|^{\frac 4{n-2}}u   \inn\ \Omega \times (0, T),  \\
u(\cdot,0) &  =u_0 \inn \Omega
\end{aligned}\right.
\label{1}\ee
In what follows we let  $n=5$. Let us fix arbitrary points  $q_1, q_2,\ldots q_k \in \Omega$.
We consider a smooth function $Z_0^*\in L^\infty (\Omega)$ with the property
 that
 $$Z^*_0(q_j) >0 \foral j=1,\ldots, k . $$

\begin{theorem} \label{teo1} Let $n=5$. For each $T>0$ sufficiently small there exists
an initial condition $u_0$ such that the solution of Problem \equ{1} blows-up at time $T$ exactly at the $k$ points $q_1,\ldots, q_k$.
It looks at main order as
$$
u(x,t) =  \sum_{j=1}^k   U_{\la_j(t), \xi_j(t)} (x)      - Z_0^*(x)  + \theta(x,t)
$$
where
$$
\la_j(t)\to 0, \quad \xi_j(t)\to q_j \ass t\to T, $$
and
 $\|\theta \|_{L^\infty} \le  T^a $ for some $a>0$.
More precisely, for numbers $\beta_j>0$ we have
$$
\la_j(t)\ =\  \beta_j (T-t)^2\, (1+ o(1))  ,  
$$
\end{theorem}

We observe that in particular, the solution predicted  by the above result has type II blow up since
 $\|u(\cdot, t) \|_{L^\infty(\R^5)} \sim   (T-t)^{-3}  $ and $\frac 1{p_S -1} = \frac 34 < 3 $.

\medskip
Our result is connected with that in \cite{schweyer} for dimension $n=4$ obtained by a very different method and only in the radially symmetric case with one blow-up point. The example in this paper is the first one  where bubbling phenomenon is obtained with arbitrary blow-up points. Moreover, the fact that $p_S=3$ for $n=4$, an odd integer, seems important in that approach.
Our computations suggest that no blow-up of this type should be present in dimensions 6 or higher.

 \medskip
 We should point out that blow-up by bubbling  (at main order time dependent, energy invariant, asymptotically singular scalings of steady states) is a phenomenon that arises in various problems of parabolic and dispersive nature.  It has for been in particular
 widely studied for the
energy critical wave equation
$$ u_{tt}  = \Delta u + |u|^{\frac 4{n-2}}.  $$
Among other works,  we refer the reader for instance to \cite{dkm,j1,kenigmerle,kst}. The method of this paper substantially differs from those in most of the above mentioned references for the parabolic case. It is close in spirit to the analysis in the works
\cite{CDM,dd,ddw,dmw,dmw1}, where the inner-outer gluing method is employed. That approach consists of reducing the original problem to solving a basically uncoupled system, which depends in subtle ways on the parameter choices (which are governed by relatively simple ODE systems).

\medskip
The rest of this paper will be devoted to the proof of Theorem \ref{teo1}.

 \section{Basics in the construction}

For notational simplicity we shall only carry out the proof in the case $k=1$ and at the end we will formulate the (relatively minor) changes needed for the general case.

\medskip
Thus we fix  a point $q\in \Omega$. Let us consider a function $Z^*_0$ smooth in $\bar \Omega$ with $Z^*_0= 0$ on $\pp \Omega$.
We assume in addition that
\be\label{negativo}
Z^*_0(q) <0. \ee
 We let $Z^*(x,t) $ be the unique solution of the initial-boundary value problem
 \be\left\{
\begin{aligned}
Z^*_t   & = \Delta Z^*    \inn\ \Omega \times (0, \infty),  \\
Z^*  & = 0 \onn \pp \Omega \times (0,\infty),\\
Z^*(\cdot,0) &  = Z^*_0 \inn \Omega.
\end{aligned}\right.
\label{11}\ee
We consider functions  $\xi(t)\to q $, and parameters $\la(t)\to 0$ as $t\to T$.
We  look for a solution of the form  \be\label{formasol} u(x,t) =  U_{\mu(t),\xi(t)}(x)  +   Z^*(x,t)  +\vp(x,t)\ee with a remainder $\vp$ which we write in the  form
\be \label{remainder}
\vp(x,t)  =   \mu^{-\frac{n-2}2} \phi\left (y,t\right ) \eta_R(y)  +  \psi(x,t)    , \quad  y = \frac{x-\xi(t)} {\mu(t)}
\ee
where
$$
\eta_R(y) =  \eta_0\left (  \frac {|y |}{ R}  \right )
$$
and $\eta_0(s)$ is a smooth cut-off function with  $\eta_0(s) =1 $ for $s<1$ and $=0$ for $s>1$.

\bigskip
Let us set
$$S(u) =  -u_t + \Delta u +  u^p  $$
and compute $S(U_{\la,\xi} + Z^*+\vp )$.
We  have that
$$
- \pp_t U_{\mu,\xi} +  \Delta U_{\mu,\xi}   =   \mu^{-\frac {n+2}2}E   -U_{\mu,\xi}^p
$$
where
\be \begin{aligned}
  E(y,t)\  = &\  \,  \mu \dot \mu [ y\cdot  \nn U(y)  + \frac{n-2}2 U(y) ] +   \mu \dot\xi \cdot \nn U (y).
\end{aligned}\label{E}\ee
Then,
setting
$$
N_{\mu,\xi} (Z) \ = \  |U_{\mu,\xi}+ Z|^{p-1}(U_{\mu,\xi}+Z)   -  U_{\mu,\xi}^p  -p U_{\mu,\xi}^{p-1}Z,
$$
 we find
$$ \begin{aligned}
S(U_{\la,\xi} + Z^* +  \vp )\ = &\     -\vp_t  + \Delta \vp    +  pU_{\mu,\xi}^{p-1}(\vp+ Z^*)  +  \mu^{-\frac {n+2}2}E  + N ( Z^*+ \vp)
\\
= &\  \  \eta_R\mu^{-\frac{n+2}2} \big [ -\mu^2 \phi_t + \Delta_y \phi   +  pU(y)^{p-1}[\phi +   \mu^{\frac{n-2}2} (Z^*+  \psi) ]  + E\, \big ]       \\
& \ -\psi_t  + \Delta_x \psi  +  p \mu^{-2}(1-\eta_R) U(y)^{p-1}( Z^*+  \psi)    + A[\phi]    \\
& \ +  B [\phi] +  \mu^{-\frac {n+2}2}E(1-\eta_R) + N (Z^* +\vp)
\end{aligned}
$$
where
$$
\begin{aligned}
A[\phi] :=  &\        \mu^{-\frac {n+2}2 } \left \{\,  \Delta_y\eta_R \phi  +      2  \nn_y\eta_R \nn_y \phi \right \}  \\
B [\phi] := &\   \mu^{-\frac{n}2 }\left \{\,  \dot\mu \big [  y\cdot \nn_y\phi  +  \frac{n-2}2\phi\big ]\eta_R  +   \dot\xi \cdot\nn_y \phi\,\eta_R  +  \big[  \dot \mu  y\cdot \nn_y \eta_R   +  \dot \xi \cdot \nn_y \eta_R\big ] \, \phi  \right  \}
\end{aligned}
$$
and we have used  $  U_{\mu,\xi}^{p-1}\vp =  \mu^{-2} U(y)^{p-1}\vp $.
Thus, we will have a solution of the full problem {\bf if }  the pair $(\phi(y,t) , \psi(x,t))$ solves the following system of equations
(which we call the outer-inner gluing system)

\be\label{I}
\mu^2 \phi_t\ = \  \Delta_y \phi   +  pU(y)^{p-1}\phi   + H( \psi, \mu,\xi)  \, \inn B_{2R} (0)\times (0,T)
\ee

\bigskip
\be\label{outer}\left \{
\begin{aligned}
\psi_t \ = &\ \Delta_x \psi  +  G(\phi,\psi, \mu,\xi)   \inn \Omega \times (0,T)\\
\psi   \ = &  -  U_{\mu,\xi} \qquad\qquad \quad \ \ \onn \pp \Omega \times (0,T) ,\\
\psi (\cdot,0) \ = &\ 0  \qquad\qquad \qquad\qquad \inn \Omega
\end{aligned}\right.
\ee
where
\be\label{HG}
\begin{aligned}
H(\psi, \mu,\xi)(y,t)\ := & \   \mu^{\frac{n-2}2} pU(y)^{p-1}   (Z^*(\xi + \mu y ,t)+  \psi(\xi + \mu y,t))  + E(y,t),
\\
G(\phi,\psi, \mu,\xi)(x,t)\ := &\  p \mu^{-2}(1-\eta_R) U(y)^{p-1}( Z^*+  \psi)    + A[\phi]
 \ +  B [\phi] \\ &\  +  \mu^{-\frac {n+2}2}E(1-\eta_R) + N ( Z^*+ \vp), \quad  y = \frac{x-\xi}{\mu}\ ,
 \end{aligned}
\ee
and $E(y,t)$ is given by \equ{E}.

\section{Formal derivation of $\mu$ and $\xi$}

Next we do a formal consideration that allows us to identify the parameters $\mu(t)$ and $\xi(t)$ at main order.
 Leaving aside smaller order terms, what approximates the inner problem \equ{I} is an equation of the form
\be
\begin{aligned}
\mu^2\phi_t  = \Delta_y \phi  & +  pU(y)^{p-1}\phi  +  h(y,t) \inn \R^n\times (0,T)   \\ \phi(y,t) &\to 0 \ass |y|\to \infty
\end{aligned} \label{kkk}\ee
with
\be
\begin{aligned}
h(y,t)\ =  & \ \mu\dot\mu \left ( U(y) + y\cdot \nn U(y) \right )    +  p \mu^{\frac{n-2}2}U(y)^{p-1} Z_0^*(q)  \\  & +\
 \mu\dot\xi\cdot  \nn U(y) + p \mu^{\frac{n}2}U(y)^{p-1} \nn Z_0^*(q)\cdot y.
\end{aligned}
\label{hh}\ee
The condition of decay in $y$ is natural in order to mitigate the effect of $\phi$ in the outer problem \equ{outer},
making at main order \equ{I} and \equ{outer} decoupled.

\medskip
At this point we recall that
the elliptic equation
$$\begin{aligned}
L[\phi] :=  \Delta_y \phi  & +  pU(y)^{p-1}\phi  =  g (y)  \inn \R^n \\  \phi(y) &\to 0 \ass |y|\to \infty,
 \end{aligned}$$
with $h(y)= O(|y|^{-2-a}) $ with $0<a<1$,  is solvable if and only if
$\phi = O(|y|^{-a})$ provided that $$
\int_{\R^n} g(y) Z_i(y)\, dy \  =\  0 \foral i= 1,\ldots, n+1,
$$
where
$$
Z_i(y) = \pp_i U (y) ,\quad i=1,\ldots, n , \quad Z_{n+1}(y) =  \frac {n-2}2 U(y) +  y\cdot  \nn U(y)
$$
These are in fact all bounded solutions of the equation  $L[Z]=0$.
 It seems reasonable to get an approximation to a solution of equation \equ{kkk}
(valid up to large $|y|$)  by  solving the elliptic equation
$$
\begin{aligned}
 \Delta_y \phi   +  pU(y)^{p-1}\phi & +  h(y,t)= 0  \inn \R^n\times (0,T)   \\ \phi(y,t) &\to 0 \ass |y|\to \infty,
\end{aligned} $$
which we can indeed do  under the conditions
\be\label{kkk1}
\int_{\R^n} h(y,t) Z_i(y)\, dy \  =\  0 \foral i= 1,\ldots, n+1,\quad t\in [0,T).
\ee
Let us check the meaning of \equ{kkk1}
for  $h$ in \equ{hh}.  Integrating against $Z_{n+1}(y)$
we get
$$
\int_{\R^n}  h(y,t) Z_{n+1}(y) \, dy \ = \   \mu\dot\mu(t) \int_{\R^N} Z_{n+1}^2dy     -\frac {n-2}2 \mu(t)^{\frac{n-2}2} Z^*_0(q) \int_{\R^n} U^p dy
$$

This quantity is zero if and only if for a certain explicit constant $\beta_n>0$ we have
$$
\dot \mu (t)  =  -\beta_n |Z^*_0(q) | \mu(t)^{\frac{n-4}2}, \quad \mu(T)= 0.
$$
For $n=5$ this equation leads to the solution
\be \label{mustar}\mu_*(t) =  \alpha  (T-t)^2, \quad \alpha =  \frac 14 \beta_n^2 |Z^*_0(q) |^2 . \ee
 In a similar way, the remaining $n$ relations in \equ{kkk1} lead us to
$\dot \xi(t) =    \mu(t)^{\frac {n-2}2} b $ for a certain vector $b$. Hence  $\dot \xi(t) =O (T-t)^3 $ and
$$ \xi(t)  =  q + O (T-t)^2.  $$
We see then  that if $u (x,t)\sim U_{\xi(t),\mu(t)}(x)$ then $$\|u(\cdot, t)\|_{L^\infty (\Omega)} \sim \mu(t)^{-\frac {n-2}2 } \sim (T-t)^{-3} \ass t\to T.
$$
This rate is type II since type I corresponds to the order $ O(T-t)^{-\frac 34}.  $

\medskip
To solve the actual parabolic problem \equ{kkk}, even assuming \equ{kkk1} further constraints are needed. Indeed,
let us recall that the operator
$L$ has a
 positive radially symmetric bounded eigenfunction $Z_0$
associated to the only positive eigenvalue $\la_0$ to the problem
$$
\label{eigen0}
L_0[ \phi ] = \la \phi  , \quad  \phi \in L^\infty(\R^n).
$$
It is known that $\la_0 $ is a simple eigenvalue and  that $Z_0$
 decays like
$$Z_0(y) \sim  |y|^{-\frac{n-1}2} e^{-\sqrt{\la_0 }\,  |y|}  \ass |y| \to \infty. $$
Let us write
$$p(t) = \int_{\R^n}  \phi(y,t)Z_0(y)\, dy, \quad   q(t) = \int_{\R^n}  h(y,t)Z_0(y)\, dy .$$
Then we compute
$$
\mu(t)^2 \dot p(t) -\la_0 p(t)  = q(t) .
$$
 Since $\mu(t) \sim (T-t)^{-2}$, then $p(t)$ will have exponential growth in time $p(t)\sim e^{\frac c{T-t}}$ unless
$$
p (t)  =      e^{  \int_0^t \frac {d\tau}{\mu^2(\tau)} } \int_t^T  e^{ - \int_0^s \frac {d\tau}{\mu^2(\tau)}}\,\mu(s)^{-2}\, q(s)\, ds
$$
This relation imposes a linear constraint  on  the initial value $\phi(y,0)$ to our desired solution $\phi(y,t)$ to \equ{I}.
\be\label{poto0}
\int_{\R^n}  \phi(y,0)Z_0(y)\, dy\  =\   \langle \ell,  h \rangle
\ee
where
$$
\langle \ell,  h \rangle\ :=\  \int_0^T  e^{ - \int_0^s \frac {d\tau}{\mu^2(\tau)}}\,\mu(s)^{-2}\, \int_{\R^n}  h(y,s)Z_0(y)\, dy\, ds
$$

\section{The linear inner and outer problems}

\subsection{The linear inner problem}
In order to deal with the inner problem \equ{I} we need to solve a linear problem like \equ{kkk} restricted
to a large ball $B_{2R}$ where orthogonality conditions like \equ{kkk1} are assumed and the initial condition of the solution depends on a scalar parameter which is part of the unknown, connected with constraint \equ{poto0}. We construct
 a solution $(\phi, \ell )$ which defines a linear operator of functions $h(y,t)$ defined on
 $$
 \DD_{2R} = B_{2R} \times (0,T)
 $$
  to the initial value problem
\begin{equation}\label{linear} \begin{aligned}
\mu^2 \phi_t\ &= \  \Delta_y \phi   +  pU(y)^{p-1}\phi   + h(y,t )  \, \inn \DD_{2R}\\
   \phi(y,0)\ &=\ \ell Z_0(y)  \inn B_{2R},
\end{aligned}
\end{equation}
for some constant $\ell $, under
the orthogonality conditions
\begin{equation}\label{uffa1}
\int_{B_{2R}} h(y,t)\,  Z_i (y) \, dy = 0  \foral i= 1,\ldots, Z_{n+1} ,\  t\in [0,T).
\end{equation}
We impose on the parameter function $\mu$ the following constraints, which are motivated on the discussion in the previous section:
let us write
$$
\mu_0 (t)  =  (T-t)^2
$$
For some positive constants $\alpha$ and $\beta$ (to be fixed later) we impose
$$
\alpha \mu_0 (t)\   \le\   \mu(t) \ \le\ \beta \mu_0 (t)  \foral t\in [0,T].
$$
Let us fix numbers $0<a<1$ and $\nu>0$.
We will consider functions $h$  that behave like
$$ h(y,t) \sim \frac {\mu_0 (t)^\nu }{ 1+ |y|^{2+a}} \inn \DD_{2R}$$
The formal analysis of the previous section would make us hope to find a solution to \equ{linear} such that
$$ \phi(y,t) \sim \frac {\mu_0 (t)^\nu }{ 1+ |y|^{a}} \inn \DD_{2R} .$$
We will find a solution so that a somewhat worse bound for $\phi(y,t)$ in space variable is found but coinciding with the expected behavior when $|y|\sim R$.   Let us define the following norms. We let  $\|h\|_{2+a,\nu}$ be the least number $K$ such that
\be
  |h(y,t) |  \ \le    \   K \frac{\mu_0 (t)^{\nu}} { 1+ |y|^{2+a}} \inn \DD_{2R}
\label{norma**}\ee
We also define  $\|\phi\|_{* a, \nu}$ be the least number $K$ with
\be
  |\phi(y,t) |  \ \le    \  K \mu_0 (t)^{\nu}  \frac {R^{n+1-a}} { 1+ |y|^{n+1}} \inn \DD_{2R}.
\label{norma*}\ee
We observe that  $\|\phi\|_{* a, \nu}\le  \|\phi\|_{a, \nu}$.

\medskip
The following is the key linear result associated to the inner  problem.
\begin{lemma} \label{lema1} There is a $C>0$ such
For all sufficiently large $R>0$ and any  $h$  with   $\|h\|_{2+a,\nu} <+\infty$
that satisfies relations $\equ{uffa1}$
there exist linear operators $$\phi = \TT^{in}_\mu [h],\quad  \ell = \ell [h] $$ which solve Problem $\equ{linear}$ and define linear operators of $h$
with
$$
|\ell [h]|  +    \| (1+|y|)\nn_y \phi \|_{* a, \nu} \ +\ \|\phi\|_{* a, \nu} \ \le \ C\,\|h\|_{ \nu, 2+a}.
$$
\end{lemma}

\proof
Letting  $\tau =  \tau_0 +  \int_0^t \mu(s)^{-2}ds $ and expressing $\phi= \phi(y,\tau)$, Problem \equ{linear} becomes
$$ \begin{aligned}
 \phi_\tau\ &= \  \Delta_y \phi   +  pU(y)^{p-1}\phi   + h(y,\tau )  \, \inn B_{2R}\times (\tau_0 ,\infty) \\
   \phi(y,0)\ &=\ \ell Z_0(y)  \inn B_{2R}.
\end{aligned}
$$
The result  then follows from Proposition 5.5 and the gradient estimates in the proof of Proposition 7.2  in \cite{CDM}.
We remark that we also have the validity of a H\"older estimate in space time with the form
\qed

\subsection{The linear outer problem}

The outer problem \equ{outer} in linear version is actually simpler than its counterpart \equ{linear}, corresponding just to the standard heat equation with
nearly singular right hand sides and zero initial and boundary conditions. Thus we consider the problem

\be\label{outer1}\left \{
\begin{aligned}
\psi_t \ = &\ \Delta_x \psi  +  g(x,t)   \inn \Omega \times (0,T)\\
\psi   \ = &  \ 0 \onn \pp \Omega \times (0,T) ,\\
\psi (\cdot,0) \ = &\ 0 \inn \Omega.
\end{aligned}\right.
\ee
The class of right hand sides $g$ that we want to take are naturally controlled by the following norms.
Let $0<a<1$, $q\in \Omega$ and $\mu_0(t)= (T-t)^2$.
We define the norms  $\| g\|_{o*}$ and  $\| \psi\|_{o}$ to be respectively the least numbers $K_1$ and $K_2$ such that
for all $ (x,t) \in  \Omega \times [0,T),$
$$ \begin{aligned}
|g(x,t) | \ &\le \  K_1\, \left [    \frac 1{\mu_0 (t)^2 } \frac 1{1 +|y|^{2+a}}  +   1   \right ]  \\
|\psi(x,t)|  \ &\le\  K_2\, \left [    \frac 1{1 +|y|^{a}}  +   T^{\frac 32a}    \right ]  \end{aligned}\ , \qquad y = \frac {x- q} { \mu_0 (t)}$$
Then the following estimate holds.

\begin{lemma}\label{lema2}
There exists a constant $C$ such that for all sufficiently small $T>0$ and any $g$ with  $ \|g \|_{o}<+\infty$,
the unique solution $\psi =  \TT^{out} [g]$  of problem \equ{outer1} satisfies the estimate
\be\label{estima}
\|\psi\|_{o*} \  \le\ C \|g \|_{o}.
\ee
\end{lemma}
\proof
We will find a positive supersolution  $\bar\psi$ to the problem
\be\left\{
\begin{aligned}
\bar \psi_t   - \Delta_x \bar \psi  \ \ge  &\   \frac 1{\mu_0^2 } \frac 1{1 +|y|^{2+a}}  +   1      \inn \Omega \times (0,T)\\
\bar \psi   \ \ge  &  \ 0 \onn \pp \Omega \times (0,T) ,\\
\bar \psi (\cdot,0) \ \ge  &\ 0 \inn \Omega.
\end{aligned}\right. \label{super} \ee
with the property that
\be\label{i3}
\bar\psi(x,t)  \ \le \   C\, \left [    \frac 1{1 +|y|^{a}}  +   T^{\frac 32a}    \right ].
\ee
After this is done,  estimate \equ{estima} follows from standard comparison.

\medskip
Let us consider the function $\psi_1(x,t) = 2 p(y) $ where $p(y)$ is a radial solution of
$$
-\Delta_y p(y)  =  \frac 1{1 +|y|^{2+a}}\inn \R^n,
$$
with $p(y) \sim \frac 1{1+ |y|^a} $ as $|y|\to +\infty $ ($p$ corresponds to the Newtonian potential of the right hand side). Then
\be \label{i1} \pp_t \psi_1   - \Delta_x \psi_1  =  \frac 1 {\mu_0^2 } \frac 1{1 +|y|^{2+a}} +  \ttt g(x,t)
\ee
where
$$
  \ttt g(x,t) =   \frac 1 {\mu_0^2 } \frac 1{1 +|y|^{2+a}}  -  2\frac {\dot\mu_0}{\mu_0} y\cdot \nn p(y).
$$
We see that for some $\gamma >0$ independent of $T$  we have
$$
 - \ttt g(x,t) \ \le\ \begin{cases} 0  & \hbox{ if }\   |x- q| \le  \gamma \sqrt{ T-t}, \\   c (T-t)^{\frac 32a -1 }  & \hbox{ if }\   |x- q| \ge  \gamma \sqrt{ T-t}     .   \end{cases}
$$
Hence, if we set $$\psi_2(x,t) =    t + c \alpha^{-1} ( T^\alpha - (T-t)^{\alpha}), \quad \alpha = \frac 32 a$$
we get
\be \label{i2}
\pp_t \psi_2   - \Delta_x \psi_2  =   c (T-t)^{\frac 32a -1 } +1  \ge     - \ttt g(x,t)    + 1.
\ee
Adding inequalities \equ{i1} and \equ{i2} we get the validity of \equ{i3} for $\bar\psi = \psi_1+ \psi_2$. \qed

\bigskip
We remark that  since for arbitrary $T'<T$,
 $$ \|g\|_\infty \le  C\mu_0(T')^{-2}\|g\|_{o} ,$$
standard parabolic estimates yield, for given $0<\alpha<1$, in addition a uniform H\"older bound in space-time of the form
\be
[ \psi]_{\alpha, T'}  \le   C\mu_0(T')^{-2} \|g\|_{o}
\label{holder} \ee
where $[ \psi]_{\alpha, T'} $ is the least number $K$ such that
$$
 |\psi(x_1,t_1) -\psi(x_2, t_2 ) |\,  \le\,  K ( |x_1- x_2|^\alpha + |t_1-t_2|^{\frac \alpha 2}  )
$$
for all $ x,x'\in \Omega, \quad 0\le t,t' \le T' $.

\section{The proof of Theorem \ref{teo1}}

With the above preliminaries we are now ready to carry out the proof of Theorem \ref{teo1} for the case $k=1$.
 We want to find a tuple $ \vec p  = (\phi,\psi, \mu,\xi)$  so that the inner-outer gluing system \equ{I}-\equ{outer} is satisfied, so that $u$ in \equ{formasol} solves \equ{1}  and the remainder \equ{remainder} is small at all times
$t\in [0,T]$.
We will achieve this by formulating the problem as a fixed point problem for $\vec  p$ in a small region of a suitable Banach space.

\medskip
We start by setting up the inner problem.
For a function $h(y,t)$ defined in $\DD_{2R}$ we write
$$
c_j[h](t) =  \frac {\int_{B_{2R} } h(y,t)\, Z_j(y)\, dy}{\int_{B_{2R}} |Z_j(y)|^2dy }
$$
so that the function
$$
\bar h(y,t) =  h(y,t) -\sum_{j=1}^{n+1} c_j[h](t)\, Z_j(y)
$$
satisfies
$$
\int_{B_{2R}} \bar h(y,t)\, Z_j(y)\, dy \ =\ 0 \foral j=1,\ldots, n+1, \quad t\in [0,T)
$$
which makes the result of Lemma \ref{lema1} applicable
to the equation
\be\left\{ \begin{aligned}
\mu^2 \phi_t\, & = \,  \Delta_y \phi   +  pU(y)^{p-1}\phi  +    {\bar H} ( \psi,\mu,\xi)   \ \  \inn \DD_{2R} \\
& \phi(\cdot ,0)  = \ell \, Z_0\inn B_{2R}
\end{aligned} \right.
\label{Inner}\ee
where
$$
{\bar H} ( \psi,\mu,\xi) =  { H} (\psi,\mu,\xi) -\sum_{j=1}^{n+1} c_j[{ H} (\psi,\mu,\xi)]\, Z_j
$$
and ${ H} (\psi,\mu,\xi)$ is defined in \equ{HG}.
Using Lemma \ref{lema1}, we find a solution to  \equ{Inner} if the following equation is satisfied.
\be \label{Inner1}
\phi = \TT_{\mu}^{in} [\,  {\bar H} ( \psi,\mu,\xi)  \,  ]  =: \FF_1 (\phi, \psi,\mu,\xi) .
\ee
We will have the inner equation \equ{I} satisfied if in addition we have
\be \label{Inner2}
c_j[{ H} (\psi,\mu,\xi)] \ =\ 0\foral  j=1,\ldots, n+1.
\ee
In addition, we have that equation \equ{outer} is satisfied provided that
\be \label{Outer}
\psi = \TT^{out} [\,  { G} (\phi, \psi,\mu,\xi)  \,  ]  =: \FF_2 (\phi, \psi,\mu,\xi) .
\ee
where the operator ${ G} (\phi, \psi,\mu,\xi)$ is defined in \equ{HG}.
We will solve System \equ{Inner}-\equ{Inner2}-\equ{Outer} using a degree-theoretical argument.

\medskip For $\la\in [0,1]$
we define the homotopy
$$
\begin{aligned}
H_\la (\psi, \mu,\xi)(y,t)\   = &\
\mu^{\frac{n-2}2} pU(y)^{p-1} Z_0^*( q)  +  \mu \dot \mu Z_{n+1}(y)\\  &\  + \mu \sum_{j=1}^n\dot\xi_j Z_j (y)  \\
& +\, \la \, \mu^{\frac{n-2}2} pU(y)^{p-1}   (
 Z^*( \xi + \mu y ,t) -Z_0^*(q)+  \psi(\xi + \mu y,t) \, )   ,\\
\end{aligned}
$$
and consider the system of equations
\be\label{homotopia}\left \{  \begin{aligned}
\phi\ =&\  \TT_{\mu}^{in} [\,  { H_\la} ( \psi,\mu,\xi)  -\sum_{j=1}^{n+1} c_j[{ H_\la} (\psi,\mu,\xi)]\, Z_j  \,  ]  
\\
c_j[ {  H_\la} &(\psi,\mu,\xi)] \ = \ 0\foral  j=1,\ldots, n+1,\\
\psi\ = &\  \TT^{out} [\, \la { G} (\phi, \psi,\mu,\xi)  \,  ] .
\end{aligned} \right.
\ee
We observe that for $\la=1$ this problem precisely corresponds to  that we want to solve,
\equ{Inner}-\equ{Inner2}-\equ{Outer}.

\medskip
It is convenient to write
$$ \mu(t)  = \mu_*(t)  + \mu^{(1)}(t), \quad \xi(t) =  q+ \xi^{(1)}(t) \quad t\in [0,T] $$
where $\mu_*(t)$ was defined in \equ{mustar},
$$\mu_*(t) =  \alpha_*  (T-t)^2, \quad \alpha_* =  \frac 14 \beta_n^2 |Z^*_0(q) |^2 $$
and $\mu^{(1)}(T)=0$, $\xi^{(1)}(t)=0$.

\medskip
We assume that we have a solution $(\phi, \psi, \mu^{(1)}, \xi^{(1)})$ of system \equ{homotopia}
where we assume the following constraints  hold:  
\be\left\{
\begin{aligned}
|\dot \mu^{(1)}(t)| +  |\dot \xi^{(1)}(t)| \  \le &\ \delta_0 \,\\
\|\phi\|_{*a, \nu}   +\|\psi\|_\infty \ \le &\ \delta_1
  \end{aligned}\right.
\label{restriccion}  \ee
where $\delta_0, \delta_1$ are  small positive constants whose value we will later adjust (and may depend on $T$)
We will also assume that $Z^*$ is sufficiently small but fixed independently of $T$, $\|Z^*\|_\infty \ll 1$.

\medskip
The function  $\mu_* (t)$ solves the equation
\be\label{polo}
\dot\mu_*(t) \int_{\R^N} Z_{n+1}^2dy     + \mu_*(t)^{\frac{n-4}2} Z^*_0(q) \int_{\R^n} pU^{p-1}Z_{n+1} dy =  0.
\ee
The equation \be c_{n+1} ( H_\la (\psi,\mu_* +\mu_1 , \xi))(t) = 0\quad t\in [0,T)  \label{poto5}\ee
which corresponds to
$$
\begin{aligned}
0 \ = & \quad \dot\mu(t) \big (\int_{B_{2R}} Z_{n+1}^2dy \big)   \ +\  \mu(t) ^{\frac{n-4}2} Z^*_0(q) \int_{B_{2R}} pU^{p-1}Z_{n+1} dy \\ &\  +\  \la\,\mu(t) ^{\frac{n-4}2} \int_{\R^n}
 pU(y)^{p-1}   (
 Z^*( \xi(t) + \mu(t) y ,t) -Z_0^*(q)+  \psi(\xi(t) + \mu(t) y,t) \, )Z_{n+1}(y)\, dy
\end{aligned}
$$
can be written as
$$
 \quad \dot\mu(t)  + \beta  \mu(t) ^{\frac{n-4}2}    =   \mu(t) ^{\frac{n-4}2}( \delta_R +   \la \theta (\psi, \xi,\mu_1) )
$$
for a suitable number $\beta>0$, $\delta_R= O(R^{-2})$  and the operator $\theta $ satisfies, for some absolute constant $C$,
$$
|\theta (\psi, \xi,\mu_1)| \ \le \ C \, [\, T\,  +\,  \|\psi\|_\infty \, ]
$$
can be written, using \equ{polo}, in the ``linearized'' form, for a suitable $\gamma>0$,
$$
\dot \mu_1  +  \frac {\gamma }{T-t} \mu_1  =    (T-t)g_0(\psi,\mu,\xi)
$$
with
$$
|g_0(\psi, \xi,\mu^{(1)},\la)(t)| \ \le \ C\, ( \|\psi\|_\infty +  T+ R^{-2} )   .
$$
The linear problem
$$
\dot \mu  +  \frac {\gamma }{T-t} \mu  =    (T-t) g(t), \quad \mu_1(T) =0
$$
can be uniquely solved by the operator in $g$,
$$
\mu(t) =  \TT^0[ g] (t) :=   - (T-t)^{-\gamma}\int_t^T   (T-s)^{\gamma+1}  g_0(s)\, ds.
$$
It defines a linear operator on $g$ with estimates
$$
\|(T-t)^{-1}\dot \mu\|_\infty +  \|(T-t)^{-2}\mu\|_\infty \le C \|g_0\|_\infty.
$$
Equation \equ{poto5} then becomes
$$
\mu^{(1)}(t)\  =\   \TT^{(0)}[\,  g_0(\psi, \xi,\mu^{(1)},\la) \, ] (t) \foral t\in [0,T)
$$
and we get
\be\label{mu1}
\|(T-t)^{-1}\dot \mu^{(1)}\|_\infty +  \|(T-t)^{-2}\mu^{(1)}\|_\infty\ \le \  C\, ( \|\psi\|_\infty +  T + R^{-2})
\ee
Similarly, equations
$$
c_j[ {  H_\la} (\psi,\mu,\xi)] \ = \ 0\foral  j=1,\ldots, n,\\
$$can be written in vector form as

\be\label{poto66}
\xi^{(1)}(t)  =  \TT^{(1)}[  g_1(\psi,\mu_1, \xi_1)] (t)\foral t\in[0,T).
\ee
where
$$
\TT^{(1)}[g] :=    \int_t^T (T-s) g(s) \, ds.
$$
and
$$
|g_1(\psi, \xi,\mu^{(1)},\la)(t)| \ \le \ C\, ( \|\psi\|_\infty +  T )   .
$$
From equation \equ{poto66} we thus find
\be\label{xi1}
\|(T-t)^{-1}\dot \xi^{(1)}\|_\infty +  \|(T-t)^{-2}\xi^{(1)}\|_\infty\ \le \  C\, ( \|\psi\|_\infty +  T )
\ee
On the other hand, we have
$$
|H(\psi, \mu,\xi)(y,t)| \ \le \    C \frac {\mu(t)^{\frac{n-2}2}} {1 +|y|^4} (\|\psi\|_\infty+ \|Z^*\|_\infty) +  \frac {\mu\dot\mu } { 1+|y|^{n-2}}  +  \frac {\mu|\dot\xi| } { 1+|y|^{n-1}}
$$
hence for any $0<a<1$ we have
$$
|H(\psi, \mu,\xi)(y,t)| \ \le \    C \frac {\mu_0(t)^{\frac{n-2}2}} {1 +|y|^{2+a}} \big (\|\psi\|_\infty+ \|Z^*\|_\infty\big ).
$$
and from the first equation in \equ{homotopia} and Lemma \ref{lema1}
we find that
\be
\| \phi\|_{*a, \nu} \ \le    \ C \big (\|\psi\|_\infty+ \|Z^*\|_\infty\big ), \quad \nu = \frac{n-2}2 .
\label{poto77}\ee
for the norm defined in \equ{norma*}.
Next we consider the last equation in \equ{homotopia}.
We recall that
$$
\begin{aligned}
G(\phi,\psi, \mu,\xi)(x,t)\ = &\  p \mu^{-2}(1-\eta_R) U(y)^{p-1}( Z^*+  \psi)    + A[\phi]
 \ +  B [\phi] \\ &\  +  \mu^{-\frac {n+2}2}E(1-\eta_R) +   N ( Z^*+  \mu^{-\frac{n-2}2}\eta_R \phi + Z^*+ \psi ),\\
 E(y,t)\  = &\  \,  \mu \dot \mu [ y\cdot  \nn U(y)  + \frac{n-2}2 U(y) ] +   \mu \dot\xi \cdot \nn U (y),\\
A[\phi] =  &\        \mu^{-\frac {n+2}2 } \left \{\,  \Delta_y\eta_R \phi  +      2  \nn_y\eta_R \nn_y \phi \right \}  \\
B [\phi] = &\   \mu^{-\frac{n}2 }\left \{\,  \dot\mu \big [  y\cdot \nn_y\phi  +  \frac{n-2}2\phi\big ]\eta_R  +   \dot\xi \cdot\nn_y \phi\,\eta_R  +  \big[  \dot \mu  y\cdot \nn_y \eta_R   +  \dot \xi \cdot \nn_y \eta_R\big ] \, \phi  \right  \}
\end{aligned}
$$

\medskip
Let us consider for example the error terms
$$ g_1(x,t) =   \mu^{-2} (1- \eta_R) U^{p-1} (Z^*+\psi), \quad g_2(x,t)= \mu^{-\frac {n+2}2}E(1-\eta_R)    .$$ We see that
$$
|g_1 (x,t) |   \ \le       \frac 1{ R^{2-\sigma} }  \mu^{-2}   \frac C{ 1 + |y|^{2+\sigma} }  (\|Z^*\|_\infty + \|\psi\|_\infty).
 $$
and
$$
|g_2 (x,t) |  \ \le \      \frac 1{\mu^2} \big[ \frac 1{|y|^{n-2}}  \mu^{-\frac {n-2}2} ( |\mu\dot\mu|+ |\mu\dot\xi|) \ \le\   \frac 1 {R^{3-\sigma}} \mu^{-2}\frac C{1+ |y|^{2+\sigma}}.
 $$
Let us now estimate the term $A[\phi]$.
Let us choose $\sigma= \frac a2$, where $a$ is the number in the definition of $\|\phi \|_{*a,\nu}$. We have
$$
\begin{aligned}
\big |A[\phi](x,t) \big|\ \le &\       \mu^{-2}  \frac 1{R^2}  \frac 1 { 1 + R^{-2-\sigma  } |y|^{2+\sigma }}\,   \mu^{-\frac{n-2}2}\sup_{R<|y|<2R} (|\phi| +  |y||\nn \phi| ) \\
\  \le &\      \mu^{-2}    \frac {R^{-\frac a2 }} {1  +  |y|^{2+\sigma }}   \|\phi\|_{*a,\frac{n-2}2}
\end{aligned}
$$
and similarly,
$$
\begin{aligned}
\big |B[\phi](x,t)\big|\ \le &\
 C \mu^{-2}   [ \mu  \dot \mu  + \mu|\dot\xi|]   \frac { R^{n+1 -a} }{ 1+|y|^{n+1} }  \|\phi\|_{*a,\frac{n-2}2} \\
 \\
\  \le &\      C\mu^{-2}    \frac { \mu^{\frac 32} R^{n+1-a}} {1  +  |y|^{2+\sigma }}   \|\phi\|_{*a,\frac{n-2}2}.
\end{aligned}
$$
 Now for some $\sigma>0$ we have
$$  \begin{aligned} \big | N ( Z^*+  \mu^{-\frac{n-2}2}\eta_R \phi + Z^*+ \psi ) \big |   \ \le &\     C
\mu^{-2}    \frac {\mu^\sigma} {1  +  |y|^{2+\sigma }}     (\|\phi\|_{*a,\frac{n-2}2} R^{n+1-a}  +  \|Z_*\|_\infty + \|\psi\|_\infty )^2 \\ & \   +  C(\|Z_*\|_\infty + \|\psi\|_\infty )^p. \end{aligned}
$$
According to the above estimates, it follows that if $T$ is sufficiently reduced depending on the large $R$ we have fixed, we get  using Lemma \ref{lema2},
\be \label{poto88}
\|\psi\|_\infty  \ \le \ C T^{\sigma'} \|Z_*\|_\infty  + R^{-\sigma'}\|\phi\|_{*a,\frac{n-2}2}.
\ee
Combining \equ{poto77} and \equ{poto88} and then using \equ{mu1}-\equ{xi1}, we finally get
\be\label{mu11}\left\{\begin{aligned}
\|\psi\|_\infty  \ \le &\ C T^{\sigma'} \|Z_*\|_\infty  \\
\|\phi\|_{*a,\frac{n-2}2} \ \le & \ C\|Z_*\|_\infty \\
\|(T-t)^{-1}\dot \xi^{(1)}\|_\infty +  \|(T-t)^{-2}\xi^{(1)}\|_\infty \ \le & \  C\, (T^{\sigma'} (\|Z_*\|_\infty +  1) + R^{-2} )\\
\|(T-t)^{-1}\dot \mu^{(1)}\|_\infty +  \|(T-t)^{-2}\mu^{(1)}\|_\infty\ \le  &\  C\, T^{\sigma'} (\|Z_*\|_\infty +  1) \end{aligned}\right.
\ee
We write System \equ{homotopia} in the form

\be\label{homotopia1}\left \{  \begin{aligned}
\phi\ =&\  \TT_{\mu}^{in} [\,  {\bar H_\la} (  \TT^{out} [\, \la { G} (\phi, \psi,\mu,\xi)  ,\mu,\xi)  
 ]  
\\
\psi\ = &\  \TT^{out} [\, \la { G} (\phi, \psi,\mu,\xi)  \,  ] \\
\mu^{(1)}\  = &\   \TT^{(0)}[\,  \ttt g_0(\psi, \xi^{(1)},\mu^{(1)},\la) \, ]\\
\xi^{(1)}\ =  &\TT^{(1)}[ \ttt g_1(\psi,\mu^{(1)}, \xi^{(1)},\la)]
\end{aligned} \right.
\ee
Here, we can write,
$$
\begin{aligned}
\ttt g_0(\psi, \xi^{(1)},\mu^{(1)},\la) = & c_R^1 \int _{B_{2R} } { H_\la} (  \TT^{out} [\, \la { G} (\phi, \psi,\mu,\xi) ] ,\mu,\xi) Z_{n+1}(y)dy\\ \ttt g_1(\psi, \xi^{(1)},\mu^{(1)},\la) = &  c_R^2 \int _{B_{2R} } { H_\la} (  \TT^{out} [\, \la { G} (\phi, \psi,\mu,\xi) ] ,\mu,\xi) \nn U (y)dy
\end{aligned}
$$
for suitable positive constants $c_R^\ell$, $\ell=0,1$.
We fix an arbitrarily small  $\ve>0$  and consider the problem defined only up to time $t= T-\ve$.
\be\label{homotopia11}\left \{  \begin{aligned}
\phi\ =&\  \TT_{\mu}^{in} [\,  {\bar H_\la} (  \TT^{out} [\, \la { G} (\phi, \psi,\mu,\xi)  ,\mu,\xi)] , \quad    (y,t)\in \bar B_{2R} \times [0,T-\ve]  
\\
\psi\ = &\  \TT^{out} [\, \la { G} (\phi, \psi,\mu,\xi)  \,  ], \quad    (x,t)\in \bar\Omega \times [0,T-\ve] \\
\mu^{(1)}\  = &\   \TT_\ve^{(0)}[\,  \ttt g_0(\psi, \xi^{(1)},\mu^{(1)},\la) \, ], \quad    t\in  [0,T-\ve] \\
\xi^{(1)}\ =  &\TT^{(1)}_\ve [ \ttt g_1(\psi,\mu^{(1)}, \xi^{(1)},\la)], \quad    t\in  [0,T-\ve] .
\end{aligned} \right.
\ee
where
$$
\TT^0_\ve [ g] (t) :=   - (T-t)^{-\gamma}\int_t^{T-\ve}   (T-s)^{\gamma+1}  g_0(s)\, ds, \quad \TT^{(1)}_\ve [g] :=    \int_t^{T-\ve} (T-s) g(s) \, ds.
$$
The key is that
the operators in the right hand side of \equ{homotopia11}  are compact when we regard them as defined in the space
of functions
$$
(\phi,\psi, \mu^{(1)}, \xi^{(1)} ) \in X_1 \times X_2\times X_3 \times X_4
$$
defined as, with their respective norms,
$$\begin{aligned}
X^1 = & \{\phi \ /\ \phi \in  C( B_{2R}\times [0,T-\ve] ), \  \nn_y \phi \in  C( B_{2R}\times [0,T-\ve] ) \}, \quad \| \phi \|_{X_1} = \| \phi \|_{\infty}+ \| \nn_y \phi \|_{\infty} \\
X^2 =  &\{\psi \ /\ \phi \in  C(\bar\Omega \times [0,T-ve] )\} ,\quad  \|\psi \|_{X_2}  =\|\psi \|_\infty  \\
X^3 = & \{ \mu^{(1)} \ /\  \mu^{(1)} \in  C^1[0,T-\ve ]  \}, \quad  \|\mu^{(1)}   \|_{X_3}=\|\mu^{(1)} \|_\infty + \|\dot \mu^{(1)} \|_\infty  \\
X^4 =  &\{ \xi^{(1)} \ /\  \xi^{(1)} \in  C^1[0,T-\ve]  \},  \quad \|\xi^{(1)} \|_{X_4}=\|\xi^{(1)}  \|_\infty+ \|\dot \xi^{(1)} \|_\infty.
\end{aligned}
$$
Compactness on bounded sets of all the operators involved in the above expression is a direct consequence of the H\"older estimate \equ{holder} for the
operator  $\TT^{out}$ and Arzela-Ascoli's theorem. On the other hand, the a priori estimate we obtained for $\ve=0$ holds equally well, uniformly  on arbitrary small $\ve>0$, and for a solution of \equ{homotopia11}.

\medskip
Leray Schauder degree applies in a suitable ball $\mathcal B $ that contains the origin in this space: essentially one slightly bigger than that defined by relations \equ{mu11}, which amounts to a choice of the parameters $\delta_0$ and $\delta_1$ in \equ{restriccion}: in fact the homotopy connects with the identity at $\la=0$, and hence the total degree in the region defined by relations \equ{mu11} is equal to 1. The existence of a solution to the approximate problem
satisfying bounds \equ{mu11} then follows. Finally, a standard diagonal argument yields a solution to the original problem with the desired size. The proof of the theorem for the case $k=1$ is concluded.

\medskip
The general case of $k$ distinct points $q_1, \ldots, q_k$ is actually identical:  in that case we have $k$ inner problems and one outer problem with analogous properties. We look for a solution of the form
$$
u(x,t) =   \sum_{j=1}^k  U_{\mu_j,\xi_j} (x) + Z^*(x,t) +    \mu_j^{-\frac{n-2}2}  \phi(y_j,t)  \eta_R(y_j) +  \psi(x,t), \quad y_j =\frac{x-\xi_j}{\mu_j}
$$
where $Z^*$ solves heat equation
with initial condition $Z_0^*$ which is chosen so that \equ{negativo} holds at all concentration points: $Z_0^*(q_j)<0$, and
$\xi_j(T)=q_j, \ \mu_j(T)=0$.

\medskip
 A string of  fixed point problems (with essentially decoupled equations associated at each point) then appears and it
solved in the same way. We omit the details.
\qed

\bigskip\noindent
{\bf Acknowledgements:}
 M.~del Pino and M. Musso have been  partly supported by grants
 Fondecyt  1160135,  1150066, Fondo Basal CMM. The  research  of J.~Wei is partially supported by NSERC of Canada.

\end{document}